\documentclass[12pt, psamsfonts, leqno]{amsart}

\usepackage{amssymb,amsfonts,amsmath,graphicx,lineno}
\usepackage{mathrsfs}
\usepackage{anysize}
\usepackage[colorlinks, citecolor=black, urlcolor=black, linkcolor=black]{hyperref}
\usepackage{url}
\usepackage{color}
\usepackage{tikz}
\usetikzlibrary{shapes, arrows}
\usepackage{pgfplots}
\tikzset{>= angle 60}

\usepackage{graphics}
\usepackage{amsmath}
\usepackage{amsthm}
\usepackage{amssymb}
\usepackage{amscd}
\usepackage{amsfonts}
\usepackage{amsbsy}
\usepackage{epsfig}
\usepackage{color}

\newtheorem{thm}{Theorem}[section]
\newtheorem{cor}[thm]{Corollary}
\newtheorem{lemma}[thm]{Lemma}
\newtheorem{prop}[thm]{Proposition}
\newtheorem{prob}[thm]{Problem}
\newtheorem{remark}[thm]{Remark}
\newtheorem{notation}[thm]{Notation}
\newtheorem{example}[thm]{Example}
\newtheorem{remarks}[thm]{Remarks}
\newtheorem{definition}[thm]{Definition}

\def\R{{\mathbb R}}

\def\N{{\mathbb N}}
\def\Z{{\mathbb Z}}

\def\cal{\mathcal }
\def\GL{{\rm GL}}
\def\R{\mathbb R} \def\Z{\mathbb Z}  

\def\N{\mathbb N}

\def\bv{{\bf v}}
\def\bx{{\bf x}}
\def\bz{{\bf z}}
\def\by{{\bf y}}

\def\norm {|\!|}

\def\norm {|\! | }

\def\<{\,<\!}
\def\>{\!>\,}

\def\SL{{\rm SL}}
\usepackage{graphics}
\usepackage{epsf}
\usepackage{epsfig}

\begin{document}

\title{Diophantine approximation with nonsingular integral transformations} 

\author{S.G. Dani and Arnaldo Nogueira}

\date{}

\maketitle

\footnote{\footnote \rm 2010 {\it Mathematics Subject Classification:}   11J20,   11J82. }





\noindent
{\bf Abstract.} Let $\Gamma$ be the multiplicative semigroup of all $n\times n$ matrices   with  integral  entries and positive determinant. 
Let $1\leq p \leq n-1$ and $V=\R^n\oplus \cdots \oplus \R^n$ ($p$ copies). 
We consider the componentwise action of $\Gamma$ on $V$. Let $\bx\in V$ be such that $\Gamma \bx$ is dense in $V$.  We discuss the effectiveness of the approximation of any target point $\by \in V$ by the orbit $\{ \gamma \bx \mid \gamma \in \Gamma\}$, in terms of $\norm \gamma \norm$, and prove in particular that for all $\bx$ in the complement of a specific null set described in terms of a certain Diophantine condition, the exponent of approximation is  $(n-p)/p$; that is, for any $\rho<(n-p)/p$, $\norm \gamma \bx - \by \norm < \norm \gamma \norm^{-\rho}$ for infinitely many $\gamma$.

\section{Introduction}

 Let $\mathcal{M}(n,\R)$, $n\geq 2$, denote the algebra of all $n\times n$ matrices $(a_{ij})$
with entries $a_{ij}$ in $\R$, and 
 $\Gamma$ be the multiplicative semigroup of all matrices in  $\mathcal{M}(n,\R)$  with  integral  entries and positive determinant.  Let $1\leq p \leq n-1$ and $\R^{(n,p)}=\R^n\oplus \cdots \oplus \R^n$ ($p$~copies), equipped with the Cartesian product topology. Consider the action of $\Gamma$ on $V$, given by the natural action on each component, by matrix multiplication on the left. 
Then for $\bx =(x_1, \dots, x_p) \in \R^{(n,p)}$, the $\Gamma$-orbit is dense in $\R^{(n,p)}$ if and only if there exists no linear combination 
$\sum_{j=1}^p \lambda_jx_j$, where $\lambda_j\in \R$ for all $j$ and $\lambda_j\neq 0$ for some $j$, which is a rational vector in $\R^n$; 
in fact the assertion holds also for the orbit of the subgroup $\SL (n,\Z)$ that is contained in $\Gamma$ (see \cite{DR}; also \cite{JSD} for the case $p=1$), and is implied by it.

When $\bx$ is such that the $\Gamma$-orbit is dense, given $\by\in \R^{(n,p)}$  and $\epsilon>0$ one may ask for $\gamma \in \Gamma$ such that $\norm \gamma \bx -\by \norm<\epsilon$, with a bound on $\norm \gamma \norm $ in terms of $\epsilon$.  
  There has been considerable interest in the 
  literature in effective results of this kind, for various group actions. In particular it was shown in \cite{LN}, for $n=2$, that given an irrational vector $\bx$ in $\R^2$ and any target vector $\by\in \R^2$  there exist a constant $C=C(\bx, \by)$ and infinitely many $\gamma $ in $\SL(2,\Z)$ such that $\norm \gamma \bx-\by \norm\leq C \norm  \gamma \norm ^{-\frac 13}$; there are also stronger results proved in \cite{LN} under some restrictions on $\by$, which we shall not go into here; see also \cite{MW}, \cite{Ke}  and  \cite{GG} for analogous results for various actions; it may be mentioned that these results are broader in their framework, but weaker in terms of the exponents involved. Here we describe some results along this theme for the action of $\Gamma$ as above; for the case $n=2$ the result is stronger in import than the result recalled above for $\SL(2,\Z)$, in the sense that for almost all initial points $\bx\in \R^2$ the corresponding statement holds  for all $\rho$ less than $1$, in place of $\rho = \frac 13$ for $\SL(2,\Z)$; see also 
  Remark~\ref{comp-remark}.
  
In the sequel we denote by $\mathcal{M}(n,\Z)$ the subring of $\mathcal{M}(n,\R)$ consisting of all matrices with entries in $\Z$. 
 For any  $\bx=(x_1, \dots, x_p)\in \R^{(n,p)}$, where $1\leq p \leq n-1$, the maximum of the absolute values of the coordinate entries of $x_j$, $1\leq j\leq p$, is called the {\it norm} of 
 $\bx$ and will be denoted by $\norm  \bx\norm  $; for a matrix $\xi \in \mathcal{M}(n,\R)$, the norm $\norm  \xi \norm $ is defined to be the norm of the $n$-tuple formed by its column vectors, or equivalently the maximum of the absolute values of the entries. For any $ \xi \in \mathcal{M}(n,\R)$  
 and a $p$-tuple $\bx=(x_1,\dots, x_p) \in \R^{(n,p)}$ we denote by $ \xi \bx$ the $p$-tuple $( \xi x_1, \dots,  \xi  x_p)$.

 We prove the following: 
 
\begin{thm}\label{main:gen}
Let $1\leq p \leq n-1$ and $\bx =(x_1, \dots, x_p)\in \R^{(n,p)}$, with 
$x_1, \dots, x_p$ linearly independent vectors in $\R^n. $ 
Let $\displaystyle 0\leq \varphi < \frac{1}{np-1}$ be such that 
 $$ 
 \inf_{\omega \in {\cal M}(n,\Z)\backslash \{0\}}  \norm \omega \bx \norm^p\norm \omega \norm^{(n-p)(1+\varphi)} >0, 
 \leqno (1.1)
$$ 
and  let   $\psi =\displaystyle{\frac p{n-p}\cdot \frac {1+n(n-p)\varphi}{1-(np-1)\varphi}}$.
Then for any $\by\in \R^{(n,p)}$ and   $\epsilon \in (0,1)$ 
 there exists a $\gamma \in \Gamma$ such that  
 $$ \norm \gamma \bx -\by \norm <\epsilon \ \hbox{ \rm and }\  \norm \gamma \norm < \epsilon^{-\psi}. \leqno (1.2)$$

\end{thm}

It is easy to see that  for any $\bx=(x_1, \dots, x_p)$ for which condition~(1.1) holds   the subspace of $\R^n$ spanned by $x_1, \dots, x_p$ contains no nonzero rational vector.  

It would be instructive to understand when condition (1.1) holds, in terms of classical notions in Diophantine approximation. Towards this we introduce the following definition.

\begin{definition} {\rm Let $1\leq p\leq  n-1$ and $\bx \in \R^{(n,p)}$. We define the {\it homogeneous exponent} of 
$\bx$, denoted by $h(\bx)$, to be the infimum of $u$ for which there exists a $c>0$ such that 
$\norm \omega \bx \norm >c\norm \omega \norm^{-u}$ for all $\omega \in {\cal M}(n,\Z)\backslash \{0\}$. }
\end{definition}

We note that a $\bx=(x_1, \dots, x_p)\in \R^{(n,p)}$ with $x_1, \dots, x_p$ linearly independent, as above, can be realised, up to a permutation of the indices, as a matrix $\left(\begin{matrix}
\xi \theta \cr
\theta 
\end{matrix} \right)$, where 
$\xi $  is a real $(n-p)\times p$ matrix  and  $\theta$ a real nonsingular $p\times p$ matrix. It turns out that then the homogeneous exponent $h(\bx)$ as above coincides with the exponent of $\xi$ in the classical sense; see Proposition~\ref{corresp}.

 \begin{cor}\label{cor:main}
 Let $1\leq p \leq n-1$. Let $\bx =(x_1, \dots, x_p)\in \R^{(n,p)}$, with 
$x_1, \dots, x_p$ linearly independent vectors in $\R^n $, be such that $h(\bx)<\displaystyle{\frac {n(n-p)}{np-1}}$ and $\by \in \R^{(n,p)}$.  
Let $$\varphi_0=\displaystyle{\frac p{n-p} h(\bx)-1} \hbox{ and } \psi_0=\displaystyle{\frac p{n-p}\cdot\frac {1+n(n-p)\varphi_0}{1-(np-1)\varphi_0}}.$$
Then for any $\psi >\psi_0$ and any $\epsilon \in (0,1)$ there exists a $\gamma \in \Gamma $ such that $$\norm \gamma \bx -\by\norm <\epsilon \hbox{ \rm and  }\norm \gamma \norm <\epsilon^{-\psi}.$$  Consequently, if $\by \notin \Gamma \bx$ then for all $\rho <1/\psi_0$ there exist infinitely many $\gamma \in \Gamma$ such that $\norm \gamma \bx -\by \norm < \norm \gamma \norm^{-\rho}$. 
 \end{cor}

 In analogy with the classical notion of very well approximable vectors we shall say that $\bx \in \R^{(n,p)}$ is {\it projectively very well approximable} if $h(\bx)$ is greater than $(n-p)/p$; see \S~4 for details. 
   From the correspondence with the classical situation noted above, viz. from Proposition~\ref{corresp}, it follows that the set of projectively very well approximable $p$-tuples $\bx$ has Lebesgue measure $0$ in $\R^{(n,p)}$. For convenience we shall also present a direct proof of this statement (see Proposition~\ref{generic}). For the tuples that are {\it not} projectively very well approximable we have the following.  
 
\begin{cor}\label{newcor}
Let $1\leq p\le n-1$ and $\displaystyle \rho < (n-p)/p$. Then for any $\bx=(x_1, \dots, x_p)\in \R^{(n,p)}$ such that $x_1, \dots, x_p$ are linearly independent and $\bx$ is not projectively very well approximable, and thus for almost all $\bx$, the following holds: for any $\rho <(n-p)/p$ and  $\by\notin \Gamma \bx$ there  exist infinitely many $\gamma \in \Gamma$  such that  
$$
\norm  \gamma \bx-\by\norm  <  \norm  \gamma \norm ^{-\rho}.
$$ 
\end{cor}

Corollary~\ref{newcor} means, in common parlance (see \S\,5 for details), that for $\bx, \by $ as in the Corollary the  exponent of approximation of the action associated to the pair $(\bx, \by)$ is at least $(n-p)/p$. We shall also show that 

\begin{thm}\label{thm:exp}
For almost all $\bx, \by$ in $\R^{(n,p)}$ the exponent is exactly $(n-p)/p$.
\end{thm}

The paper is organized as follows. In the next section we prove a result on intersections of affine lattices with certain special sets being nonempty, on which the proof of the main theorem is based. Theorem~\ref{main:gen} is  proved in \S\,3.  In \S\,4 we discuss the relation between the homogeneous exponent and the classical exponents, and related issues of approximability, and prove  Corollaries~\ref{cor:main} and~\ref{newcor}. Theorem~\ref{thm:exp} is proved in \S\,5.

\section{A result on affine lattices in $\R^d$}

Towards the proof of Theorem~\ref{main:gen} we first prove in this section a result on intersection of affine lattices in $\R^d$ with parallelopipeds,  Proposition~\ref{genprop}. The proof of the proposition is by application of Theorem IV of \cite{meyer}. Here we consider $\R^d$ as a $d$-dimensional vector space over $\R$, with a fixed basis $\{e_1,\dots , e_d\}$.  We denote by  $\Z^d$ the lattice consisting of integral vectors with 
respect to the basis $\{e_1, \dots , e_d\}$.

\begin{prop}\label{genprop}
Let $V=\R^d$, with $d\geq 3$, and let $V_1$ and $V_2$ be vector subspaces of $V$ of dimensions $d_1\geq 2$ and $d_2\geq 1$ such that $V=V_1\oplus V_2.$ Suppose  that there exist $\displaystyle \delta \in \left(0, \frac{d_2}{d_1-1}\right)$ and $0<\kappa \leq 1 $ 
such that   for any $z=u+w\in \Z^d\backslash
\{0\}$, with $u\in V_1$ and $w\in V_2$,
$$
\norm  u\norm ^{d_1} \, \norm  z\norm^{d_2+\delta}>\kappa,
 \leqno (2.1) 
 $$
and let $\chi=d_1(1+\delta)/(d_2 -d_1\delta +\delta)$.
Let $R_1$ and $R_2$ be  compact  
convex subsets of $V_1$ and $V_2$ respectively, with nonempty interiors in the respective subspaces, and  for all $s,t>0$ let  
 $$
 \Omega (s,t) = \{v=u+w\in \R^d \mid  u\in sR_1, w\in tR_2 \}.
 $$ 
Then there exist  constants $\sigma >0$ and $\epsilon_0>0$  such that for all $0<\epsilon <\epsilon_0 $ and all $v\in \R^d$, $\Omega (\epsilon, 
\sigma \epsilon^{-\chi})\cap (v+\Z^d)\neq   \emptyset.$ 
\end{prop}

\proof  The statement is independent of the norm, and hence by modifying the norm, for convenience, we may assume that for any $u\in V_1$ and $v\in V_2$ we have $\norm u+v\norm =\max \{\norm u\norm, \norm v\norm\}$, and that 
 $R_1$ and $R_2$ are contained in $B (0,\frac 12)$, the open ball in $\R^d$ with radius $\frac 12$ and center at $0$.  
 
Let $\ell $ be the Lebesgue measure on $V$ such that $\{\sum_{j=1}^dt_je_j\mid  t_j\in [0,1] \hbox{ \rm for all } j \}$ has measure $1$. We note that if $E$ is a compact subset  such that the set of differences $E-E:=\{x-y\mid x,y\in E\}$ contains no nonzero point of $\Z^d$ then $\ell (E)< 1$. 
 
 Now let $0<a< \kappa $ be arbitrary and $S=\Omega (a,\kappa a^{-d_1/(d_2+\delta)})$ and  $S'=S-S$. Consider any $y=u+w\in S'$, with $u\in V_1$ and $w\in V_2$. If $\norm y \norm = \norm w \norm$ then 
 $$
 \norm u\norm ^{d_1} \norm y\norm ^{d_2+\delta}=\norm u\norm ^{d_1} \norm w\norm ^{d_2+\delta}< a^{d_1}(\kappa a^{-d_1/(d_2+\delta)})^{d_2+\delta}=
\kappa^{d_2+\delta}\leq \kappa, 
$$ while on the other hand if $\norm y\norm =\norm u\norm$ then we have $\norm u\norm ^{d_1}\norm y\norm^{d_2+\delta}= \norm u\norm^{d_1+d_2+\delta}<\kappa$. 
Hence by the condition in the hypothesis $S'$ does not contain any nonzero element of
$\Z^d$.  Since $S$ is a compact subset, by the observation above this implies that $\ell(S)<1$. 
 
Let $m=[\ell(S)^{-1}]+1$, the smallest integer exceeding $\ell (S)^{-1}$. Then by \cite{meyer}, Theorem~IV, page 9,  $mS \cap (v+\Z^d)\neq \emptyset$ for all $v\in \R^d$.  We shall deduce from this the desired assertion as in the Proposition. 
 
Let $ l_1$ and $l_2$  denote the Lebesgue measures on  $V_1$ and $V_2$ respectively such that $l_1(R_1)=l_2(R_2)=1$. There  exists $\lambda>0$ such that $\ell=\lambda (l_1\times l_2)$.  Then we have   
$$
\ell(S)=\lambda a^{d_1}\cdot (\kappa a^{-d_1/(d_2+\delta)})^{d_2}=\theta a^{d_1\delta/(d_2+\delta)},
$$ 
where $\theta =  \lambda  \kappa^{d_2}$.  As $m=\ell(S)^{-1}$ and $l(S)<1$, we have  $m<2(\theta a^{d_1\delta/(d_2+\delta)})^{-1}$. It follows that the set $mS$, which equals $\Omega (ma,m\kappa a^{-d_1/(d_2+\delta)})$, is contained in the set $$E_a:=\Omega (2\theta^{-1}a^{1- \frac{d_1\delta}{d_2+\delta}},2\theta^{-1}\kappa a^{-\frac {d_1(1+\delta)}{d_2+\delta}}), $$ 
and hence $E_a$ also intersects $v+\Z^d$ nontrivially for all $v\in \R^d$, for each $a\in (0,\kappa)$. 

We now show that the desired assertion holds for the choices 
$$\sigma =2 \kappa \theta^{-(1+\chi)} \hbox{ \rm and }\epsilon_0= \theta^{-1}\kappa^{1- \frac{d_1\delta}{d_2+\delta}};$$  
to that end we  prove  that for any $\epsilon \in 
(0,\epsilon_0)$ there exists $a\in (0,\kappa)$ such that the set $E_a$ as above  is contained in $\Omega (\epsilon, \sigma \epsilon^{-\chi})$, which by the preceding observation yields the desired conclusion. 
Let 
$\epsilon \in (0, \epsilon_0)$ be given. Since $d_1>1$ and $\delta <\frac {d_2}{d_1-1}$, $d_1\delta < d_2+\delta$ and hence there exists  $0<a<\kappa $ such that $\theta^{-1}a^{1- \frac{d_1\delta}{d_2+\delta}}=\epsilon$. For this choice of $a$ we have  
$$2\theta^{-1}\kappa a^{-\frac {d_1(1+\delta)}{d_2+\delta}}=
\sigma \theta^\chi a^{-\frac {d_1(1+\delta)}{d_2+\delta}}
=\sigma \epsilon^{-\chi}a^{(1- \frac{d_1\delta}{d_2+\delta})\chi}a^{-\frac {d_1(1+\delta)}{d_2+\delta}}
= \sigma \epsilon^{-\chi},$$ as $\chi =\frac {d_1(1+\delta)}{(d_2-d_1\delta +\delta)}$. Applying the observation above for this $a$ we get that the corresponding set $E_a$ is contained in $\Omega (\epsilon, 
\sigma \epsilon^{-\chi})$ and consequently $\Omega (\epsilon, \sigma\epsilon^{-\chi})\cap (v+\Z^d)$ is nonempty for all $v\in \R^d$.   This proves the proposition. \qed

\section{Proof of Theorem~\ref{main:gen}}

The proof will be by application of the Proposition~\ref{genprop} to the vector space 
$V=\mathcal{M}(n,\R)$, realized as $\R^d$ with $d=n^2$, and $\Z^d$ identified with $\mathcal{M}(n,\Z)$. 
We follow  the notation as in the statement of the theorem. 
Let $x_1, \dots, x_p \in \R^n$ be as in the hypothesis and let $x_{p+1}, \dots, x_n\in \R^n$ be chosen so that $x_1, \dots, x_n$ are linearly independent.  

For each $i,j \in \{1, \dots, n\}$ let $\beta_{ij}\in {\cal M}(n,\R)$ be the matrix such that for all $k\in \{1, \dots, n\}$, $\beta_{ij}x_k=x_i$ if $k=j$ and $0$ otherwise. For each $j=1,\dots , n$ let $S_j$ be the subspace of ${\cal M}(n,\R)$ spanned by 
$\{\beta_{ij}\mid i=1, \dots, n\}$.
Let $V_1=\sum_{j=1}^pS_j$ and $V_2=\sum_{j=p+1}^n S_j$. Then $V_1$ and $V_2$ are vector subspaces, and as $x_1, \dots, x_n$ are linearly independent it follows that $V_1$ and $V_2$ are of dimensions $d_1=np$ and $d_2=n(n-p)$ respectively and $V=V_1\oplus V_2$. 
On  $V$ we define a (new) norm $\norm \cdot \norm_V$ by setting 
$$\norm \xi \norm_V= \max_{1\leq j \leq n} \norm \xi x_j\norm, \  \hbox { \rm for all } \xi \in {\cal M}(n,\R). $$ 
By linear independence of  $x_1, \dots, x_n$ there exists a $c\geq 1$ such that for all $\xi\in {\cal M}(n,\R)$, 
$$ c^{-1}\norm \xi \norm \leq \norm \xi \norm_V \leq c\norm \xi\norm .\leqno (3.1) $$ 
We note also that for $\xi={\bf v}_1+{\bf v}_2$ with ${\bf v}_1\in V_1$ and ${\bf v}_2\in V_2$, we have  
$$\norm {\bf v}_1 \norm_V =\max_{1\leq j \leq n} \norm {\bf v}_1x_j\norm =\max_{1\leq j \leq p} \norm \xi x_j\norm=\norm \xi \bx \norm. \leqno (3.2)$$ 

Now let $\varphi$ be as in the hypothesis of the theorem and let $\delta=n(n-p)\varphi$.  Then $\delta \in (0, d_2/(d_1-1))$.
 By condition~(1.1) there exists $\kappa_1 >0$ such that
 $$
  \norm \omega \bx\norm^{np} \norm \omega\norm^{n(n-p)(1+\varphi)}>\kappa_1 
 \hbox{ \rm for all } \omega \in {\cal M}(n,\Z)\backslash \{0\}.  
 \leqno (3.3)$$
 We recall that  $np=d_1$ and $n(n-p)(1+\varphi)=n(n-p)+\delta =d_2+\delta$. In view of (3.1) and (3.2), (3.3) therefore implies that there exists a constant $\kappa >0$ such that for $\omega \in {\cal M}(n,\Z)\backslash \{0\}$, if $\omega ={\bf v}_1+{\bf v}_2$, with ${\bf v}_1\in V_1$ and ${\bf v}_2\in V_2$, then  
 $$
  \norm {\bf v}_1\norm_V^{d_1} \norm \omega\norm_V^{d_2+\delta}>\kappa 
 \hbox{ \rm for all } \omega \in {\cal M}(n,\Z)\backslash \{0\}.  
$$
Hence condition~(2.1)
  in the hypothesis of Proposition~\ref{genprop} is satisfied for  $V_1$,  $V_2$ and $\delta$ as above. We note that in this case $\chi$ as in the Proposition is given by 
  $$\chi=\frac{d_1(1+\delta)}{d_2-d_1\delta +\delta}=\frac{np(1+n(n-p)\varphi)}{n(n-p)(1-(np-1)\varphi))}=\frac p{n-p}\cdot \frac {1+n(n-p)\varphi}{1-(np-1)\varphi}=\psi,$$
with  the last term $\psi$  as defined in the statement of the theorem.  
We shall apply the conclusion of the Proposition in this case for the  choices of compact subsets as described below.

Now let $\by =(y_1, \dots, y_p)$, $y_j\in \R^n$, $j=1,\dots, p$, be given.  Let $q$ be the rank of $(y_1, \dots, y_p)$, namely the maximal number of linearly independent $y_j$'s; by re-indexing we shall assume, as we may, that $y_1, \dots, y_q$ are linearly independent. 

We shall now first consider the case with $y_j=0$ for all $j=q+1, \dots, p$. 
Let $U$ and $W$ be the subspaces defined by 
$$
U=\sum_{j=1}^q S_j \hbox{ \rm and } W=\sum_{j=q+1}^p S_j;
$$ 
we note that $V_1= U + W$.

Now let  $g_0\in {\cal M}(n,\R)$ be the (unique) element such that 
$g_0x_j=y_j$ for all $j=1,\dots , p$ and  $g_0x_j=0$ for $j=p+1, \dots, n$. Then $g_0\in V_1$. Let 
$g_0=g_1+g_2$ be its decomposition with $g_1\in U$ and $g_2\in W$. Let 
$$\Theta =\left \{ \sum_{j=1}^n u_j \in {\cal M}(n,\R)\mid u_j\in S_j \hbox{ \rm and }\norm u_j \norm < \frac 1n\right \}.$$ 
Since by assumption 
$y_1, \dots, y_q$ are linearly independent,  $g_1$ has rank $q$. We can choose $\eta \in W\cap \Theta$ with  rank $n-q$, so that  $\det (g_1+\eta)\neq 0$, and by adjusting the sign in one of the columns of $\eta$ we can further arrange so that $\det (g_1+\eta)>0$. Using the continuity of the determinant function we conclude that there exist  neighbourhoods $N$ and $K$ of $0$ in $U$ and $W$ respectively,  such that $\det (g_1+\phi + \eta +\psi)> 0$ for all $\phi \in N$ and $\psi \in K$; we shall further choose $N$ and $K$ to be compact and convex, contained in $\Theta$, and such that $-\eta \notin K$; we note that since the rank of $\eta$ is $n-q$, in particular it is a non-zero element. 

Let $\eta =\sum_{j=q+1}^n \eta_j$, with $\eta_j\in S_j$, be the decomposition of $\eta$ as above. For each $j=1,\dots , n$ let $I_j$ be a compact convex subset of $S_j$ satisfying the following conditions:

i) for $j=1, \dots, q$, $I_j$ is a compact neighbourhood of $0$ in $S_j$, contained in $\frac 1n N$;

ii) if $j=q+1, \dots, n$, $I_j$ is a compact neighbourhood of $\eta_j$ in $S_j$, contained in $\eta_j+\frac 1n K$.

For application of Proposition~\ref{genprop} we now choose 
$R_1=\sum_{j=1}^pI_j$ and $R_2=\sum_{j=p+1}^n I_j$. We note that $R_1$,  and $R_2$ are  compact convex subsets of $V_1$ and $V_2$,  with nonempty interior in the respective subspaces. 
Thus the condition in the proposition is satisfied for $R_1, R_2$. 
As in Proposition~\ref{genprop}, for any positive real numbers $s, t$ let 
$$\Omega (s, t)= \{\bv= \bv_1+\bv_2\mid \bv_1\in sR_1,\bv_2\in tR_2 \}.$$ Then by the proposition  there exist constants $\sigma >0$ and $\epsilon_0>0$ such that for any $\epsilon \in (0,\epsilon_0)$ 
and $w\in \R^d$ we have $\Omega (\epsilon, \sigma \epsilon^{-\psi})\cap (w+{\cal M}(n, \Z))) \neq \emptyset. $ We shall also assume, as we may that $\sigma \geq \epsilon_0^{1+a}$.

 We choose $w=-g_1$.  Let  $\epsilon \in (0,\epsilon_0)$ be given. 
Then $\Omega (\epsilon, \sigma \epsilon^{-\psi}) \cap (-g_1+{\cal M}(n, \Z)) \neq \emptyset,$ and hence  there exist $\theta \in \Omega (\epsilon, \sigma \epsilon^{-\psi})$ and $\gamma \in {\cal M}(n, \Z)$ such that $\theta =-g_1+\gamma$. Let $\theta =\sum_{j=1}^n \theta_j$, where $\theta_j\in S_j$ be the decomposition of $\theta $ in $\R^d$. Then from the definition of the sets we get that for $\theta_j\in \epsilon I_j$ for all $j=1, \dots, p$ and 
$\theta_j\in \sigma\epsilon^{-\psi}I_j$ for $j=p+1, \dots, n$.

We now show that the inequalities (1.2) as in the theorem hold for this $\gamma$.  Consider first $1\leq j \leq p$. 
The choice of $g_1$ as the $U$-component of $g_0$, implies that $g_1x_j=y_j$ if $j=1, \dots, q$ and $g_1x_j=0$ if $j=q+1, \dots n$. Also, by assumption we have $y_j=0$ for $j=q+1, \dots, p$. Together this implies that $y_j=g_1x_j$ for all $j=1, \dots , p$.
Also, for these $j$ we have  $\theta_j\in \epsilon I_j \subset \epsilon \Theta$, and hence $\norm \theta_j \norm <\epsilon /n$. Thus
$$\norm \gamma x_j -y_j\norm =\norm \gamma x_j-g_1x_j\norm =\norm (\gamma -g_1)x_j\norm =\norm \theta x_j\norm=\norm \theta_j x_j\norm\leq 
n \norm \theta_j\norm \norm x_j\norm 
< \epsilon, \leqno (3.4)$$
as $\norm x_j\norm =1$. Now consider $p+1 \leq j \leq n$. 
Then we have $g_1x_j=0$, so $\gamma x_j=\theta x_j=\theta_jx_j$ and since $\theta_j\in \sigma\epsilon^{-\psi}I_j\subset \sigma\epsilon^{-\psi}\Theta$ we get 
$$\norm \gamma x_j\norm =\norm \theta_jx_j\norm \leq  n\sigma\epsilon^{-\psi}\norm \theta_j\norm \norm x_j\norm  <\sigma\epsilon^{-\psi}, \leqno (3.5)$$ since $\norm 
x_j\norm =1$. Since by choice $\sigma \geq \epsilon^{1+a}$, the inequalities (3.4) and (3.5) together imply also  that $\norm \gamma \norm <\sigma\epsilon^{-\psi}$. 
This shows that the inequalities (1.2)  in the statement of the theorem hold for the matrix $\gamma \in {\cal M}(n,\Z)$. 

We shall now  show that  $\gamma\in \Gamma$, namely that $\det \gamma >0$. 
Consider the element $$
\gamma'= g_1 +\sum_{j=1}^q \theta_j +\sum_{j=q+1}^p\epsilon^{-1}\theta_j +\sum_{j=p+1}^n \sigma^{-1}\epsilon^a \theta_j.
$$ 
 For $j=1, \dots q$, $\theta_j\in I_j\subset  \frac 1n N$, 
 and since $N$ is a convex neighbourhood of $0$ in $U$ it follows that 
$\sum_{j=1}^q\theta_j \in N$. For $j=q+1, \dots, p$ we have $\epsilon^{-1}\theta_j\in I_j \subset \eta_j+\frac 1n K$, and similarly for $j=p+1, \dots, n$, $\sigma^{-1}\epsilon^{a}\theta_j\in I_j \subset \eta_j+\frac 1n K$. Recalling that $K$ is a convex neighbourhood of $0$ in $ W$  we deduce from this that
$$\sum_{j=q+1}^p\epsilon^{-1}\theta_j + \sum_{j=p+1}^n \sigma^{-1}\epsilon^a \theta_j \in \sum_{j=q+1}^n\eta_j +K=\eta +K.$$  Altogether we get that $\gamma'$ is an element of the form $g_1+\phi +\eta +\psi$, with $\phi \in N$ and $\psi \in K$. By the choices of $N$ and $K$ this implies that $\det \gamma' >0$. 
We now note that $\gamma x_j=\gamma' x_j$ for $j=1,\dots ,q$,   $\gamma x_j=\epsilon\gamma'x_j$ for $j=q+1, \dots, p$ and 
$\gamma x_j=\sigma\epsilon^{-\psi}\gamma'x_j$ for $j=p+1, \dots, n$.
Since $x_1, \dots, x_n$ is a basis of $\R^n$ this implies that $\det \gamma = \epsilon^{p-q}\cdot \sigma^{n-p}\epsilon^{-(n-p)a}\det \gamma'$, showing that $\det \gamma >0$ as sought to be proved. 
This proves the theorem in the case under consideration, namely when $y_j=0$ for $j=q+1, \dots, p$. 

Now consider the general case, with $y_j$ possibly nonzero for $q+1\leq j\leq p$. Let $\by_0=(y_1,\dots, y_q, 0, \dots, 0)$, (with $p-q$ zeros inserted).
There exists a nonsingular $p\times p$ matrix $\theta$ such that    $\by=\by_0\theta$.
Let $\tilde \bx= \bx \theta^{-1}$.  It is straightforward to see that the condition in  Theorem~\ref{main:gen}  involving (1.1) holds for $\tilde \bx$ in place of $\bx$. 
Applying the special case as above to $\tilde \bx$, with $\by_0$ in place of $\by$, we get that there exists a constant $\sigma$ such that for any $\epsilon \in (0,1)$, there 
exists $\gamma \in \Gamma$ such that $\norm \gamma \tilde \bx - \by_0\norm <\epsilon $ and $\norm \gamma \norm\leq \sigma
\epsilon^{-\psi}$.   There exists a constant $\alpha\geq 1$ such that 
for any $n\times p$ matrix $\xi$, $\norm \xi \theta \norm \leq \alpha \norm \xi\norm$, and thus we get
$$
\norm \gamma \bx -\by \norm =\norm \gamma \tilde \bx\theta - \by_0\theta\norm \leq \alpha \norm \gamma \tilde \bx -\by_0\norm <\alpha \epsilon \hbox{ \rm and } \norm \gamma \norm <\sigma \epsilon^{-\psi}.
$$
Choosing such a $\gamma$ for $\epsilon /\alpha$ in place of $\epsilon$ we get $\gamma$ such that $\norm \gamma \bx -\by\norm <\epsilon $ and $\norm \gamma \norm < C\epsilon^{-\psi}$ where $C=\sigma \alpha^{\psi}$. This proves the assertion in the theorem  in the general case as well.  \qed

\bigskip

\section{Homogeneous exponents and projective approximability}

Let $1\leq p \leq n-1$ and $q=n-p$. For any natural numbers $k,l$ we denote by $\Z^{(k,l)}$ the lattice in $\R^{(k,l)}$ (notation as in \S\,1) consisting of elements whose coordinates are integers.  
We recall that for any $\xi \in \R^{(q, p)}$ the Diophantine exponent $e(\xi)$, in the classical sense, is  the supremum of all $a$ such that 
$$\inf_{\beta \in {\mathcal M} (p, \Z)} \norm \alpha \xi +\beta\norm <\norm \alpha\norm^{-a} \hbox{ \rm for infinitely many }\alpha\in \Z^{(p, q)}.$$   

 Let  $\xi \in \R^{(q, p)}$ be given. For $\alpha \in \Z^{(p, q)}$ we define 
$$d(\alpha) =\inf_{\beta \in {\mathcal M}(p,\Z)} \norm \alpha \xi +\beta \norm.$$ We note that if for some $a$, there exists $c>0$ such that $d(\alpha)>c\norm \alpha\norm^{-a}$ for all $\alpha \in 
\Z^{(p,q)}\setminus \{0\}$ then $a>e(\xi)$, and conversely if $a>e(\xi)$ then there exists $c>0$ such that $d(\alpha)>c\norm \alpha\norm^{-a}$ for all $\alpha \in \Z^{(p,q)}\setminus \{0\}$. Thus $e(\xi)$ is the infimum of $a$ such that for some $c>0$ we have $d(\alpha )>c\norm \alpha \norm^{-a} $ for all $\alpha \in \Z^{(p,q)}\setminus \{0\}$.

\begin{prop}\label{corresp}
Let $\xi \in {\mathcal M}(q\times p, \R)$,  $\theta  \in \GL (p,\R)$, and $\bx =\left(\begin{matrix}
\xi \theta \cr
\theta 
\end{matrix} \right)$. Then $h(\bx)=e(\xi)$. In particular $\bx $ is projectively very well approximable if and only if $\xi$ is very well approximable. 

\end{prop}

\proof 
It is easy to see that the homogeneous exponents of $\left(\begin{matrix}
\xi \theta \cr
\theta 
\end{matrix} \right)$ and  $\left(\begin{matrix}
\xi  \cr
I
\end{matrix} \right)$, where $I$ is the $p\times p$ identity matrix, are the same. Hence we may assume, as we shall, that $\theta =I$. 

We now write $\omega \in \mathcal{M}(n,\Z)\setminus \{0\}$ in the form $(\alpha, \beta)$, with $\alpha \in  {\cal M}(q\times p, \Z)$ and $\beta \in {\cal M}(p,\Z)$, expressed canonically. 
Let $b\geq 0$ be arbitrary. It is easy to see that 
$$ \inf_{\omega \in \mathcal {M} (n,\Z)\setminus \{0\}}\norm \omega \bx\norm \norm \omega \norm^{b} 
=\inf_{\alpha\in \Z^{(p,q)} \setminus \{0\}, \beta \in {\cal M}(p,\Z), \norm \alpha \xi +\beta \norm  \leq 1}\norm \alpha \xi +\beta \norm\norm (\alpha, \beta)\norm^{b}. 
$$
When $\norm \alpha \xi +\beta\norm \leq 1$ we have $\norm \beta \norm \leq \norm \alpha \xi\norm +1 \leq \norm \alpha \norm \norm \xi \norm +1 \leq (\norm \xi \norm +1) \norm \alpha \norm$.  Hence we get that 
$$ \inf_{\alpha \neq 0} d(\alpha) \norm \alpha \norm^{b}\leq \inf_{\omega \in \mathcal {M} (n,\Z)\setminus \{0\}}\norm \omega \bx\norm\norm \omega \norm^{b} \leq (\norm \xi \norm +1)\inf_{\alpha \neq 0} d(\alpha) \norm \alpha \norm^{b}.\leqno (4.1) $$
Then $h(\bx)$ is by definition the infimum of $b$'s for which the middle term in the above inequalities is positive, while by the observation preceding the proposition the infimum of $b$'s for which the extreme terms are positive is $e(\xi)$. Hence we get that $h(\bx)=e(\xi)$. This proves the first assertion in the Proposition. The second assertion follows  immediate from the first, since $\bx$ being projectively very well approximable is defined by the condition that $h(\bx)>q/p$, while $\xi$ being very well approximable corresponds to $e(\xi)>q/p$. \qed

It is well known that very well approximable matrices (viewed as vectors) $\xi$ in $\R^{(q, p)}$ form a set of Lebegue measure $0$ in the latter space. From the correspondence as above it follows that the set of projectively very well approximable $\bx$ form a set of $0$ Lebesgue measure in $\R^{(n, p)}$. We include here a direct proof of this for the convenience of the reader.

\begin{prop}\label{generic}
Let $1\leq p \leq n-1$. Then the set of $\bx$ in $\R^{(n, p)}$ which are projectively very well approximable has Lebesgue measure $0$ in $\R^{(n, p)}$.
\end{prop}

\proof Let $1\leq p\leq n$ and  
$S=\{\bx =(x_1, \dots, x_p)\in \R^{(n,p)}\mid  \norm \bx \norm \leq 1\}$. We denote by  $\nu$ be the restriction of the Lebesgue measure on $\R^{(n,p)}$ to $S$. Let $\chi > 0$ be given and let $S'=\{\bx \in S\mid \inf_{{\cal M}(n,\Z)\backslash \{0\}}  \norm \omega \bx \norm^p\norm \omega \norm^{n-p+\chi}=0\}$. To prove the first assertion of the Proposition clearly it 
suffices to show that $S'$ has measure $0$. 
  
For $r=1, \dots , n$ let ${\cal M}_r$ denote the set of matrices in $ {\cal M}(n,\R)$ with rank $r$, and ${\cal M}_r(\Z)$ the subset consisting of all integral matrices in it.  
 It is straightforward to verify that there exists a constant $c>0$ such that for all $\xi \in {\cal M}_r$ with   $\norm \xi \norm =1$, for any $\theta >0$ we have $$\nu (\{\bx \in S\mid \norm \xi \bx \norm <\theta\} )\leq c\theta^{rp}. \leqno (4.2)$$
For any $\omega \in {\cal M}(n,\Z)$ and $\epsilon\in (0,1)$ let 
$$
S(\epsilon, \omega)=\{\bx \in S\mid \norm \omega \bx \norm^p  \norm \omega \norm^{n-p+\chi}< \epsilon\}.
$$
Then for any $\omega \in {\cal M}_r(\Z)\backslash \{0\}$ and $\bx\in S(\epsilon, \omega)$ we have $\norm \frac {\omega}{\norm \omega \norm} \bx \norm  < \epsilon \norm \omega \norm^{-(n+\chi)/p}$, and hence by~(4.2) we get  
$$
\nu (S(\omega, \epsilon))\leq c(\epsilon\norm \omega \norm^{-(n+\chi)/p})^{rp}\leq c\epsilon \norm \omega \norm^{-(nr+\chi)}. 
\leqno (4.3) 
$$ 
We fix $1\leq r \leq n$ and for $q\in \N$ let $$
N_q=\#\{\gamma \in {\cal M}_r(n,\Z)\mid  \norm \gamma \norm = q\}, 
$$
the cardinality of $N_q$. Then it follows from the second assertion in  Theorem 1 of~\cite{Ka} that  there exists a positive constant constant $C=C(n,r)$ such that, for every $q\in \N$,
$$
N_1+\ldots + N_{q-1} \le Cq^{nr} \log q.
\leqno{(4.4)}
$$
Together with (4.3) and (4.4) this implies that for all $r=1, \dots, n$ and $q\in \N$ we have 
$$
\sum_{k=1}^q \sum_{\omega \in {\cal M}_r(\Z), \norm \omega \norm =k}\nu (S(\epsilon, \omega))
\le \sum_{k=1}^q N_k  \frac{c\epsilon}{k^{nr+\chi}} =c\epsilon \sum_{k=1}^q  \frac{N_k }{k^{nr+\chi}}. $$
Rewriting the right hand side of the preceding inequality we  obtain
$$
 \sum_{k=1}^q   \frac{N_k}{k^{nr+\chi}}
=  \sum_{k=1}^{q-1} (N_1+\ldots + N_k) \left( \frac{1}{k^{nr+\chi}} - \frac{1}{(k+1)^{nr+\chi}} \right)+ 
\frac{N_1+\ldots + N_q}{q^{nr+\chi}} .
$$
Using the mean value Theorem, we get $\displaystyle \frac{1}{k^{nr+\chi}} - \frac{1}{(k+1)^{nr+\chi}} < \frac{1}{k^{nr+1+\chi}}$, thus
$$
\sum_{k=1}^q   \frac{N_k}{k^{nr+\chi}} \le  \sum_{k=1}^{q-1} (k+1)^{nr} \log (k+1)  \frac{1}{k^{nr+1+\chi}} + \frac{(q+1)^{nr} \log (q+1)}{q^{nr+\chi}} 
$$
$$
=  \sum_{k=1}^{q-1} \left(1+ \frac{1}{k}  \right)^{nr} \frac{ \log (k+1)}{ k^{1+\chi} } 
+      \left(1+ \frac{1}{q} \right)^{nr}     \frac{\log (q+1)}{q^{\chi}}. 
$$
Since $\chi >0$, this shows that 
$\sum_{\omega \in {\cal M}(n,\Z)\backslash \{0\}} \nu (S(\epsilon, \omega))<\infty $. Hence by the Borel-Cantelli lemma for almost all $\bx \in S$, $\bx$  is contained in $S(\epsilon, \omega)$ for at most finitely many $\omega$'s. Hence  we get that $\nu (S')=0$, as sought to be proved.  \qed

{\it Proof of Corollary~\ref{cor:main}}:

We follow the notation as in the hypothesis of the Corollary. Let $\psi >\psi_0$ be given. Then there exists $\varphi \in (\varphi_0, 1/(np-1))$ such that $\displaystyle{\psi \geq \frac p{n-p}\cdot \frac {1+n(n-p)\varphi}{1-(np-1)\varphi}}$. 
We note that $\displaystyle{\frac {n-p}p(1+\varphi) >h(\bx)}$. From the definition of the homogeneous exponent $h(\bx)$ this implies that condition~(1.1) of Theorem~\ref{main:gen} is satisfied for $\varphi$. The first statement in the corollary therefore follows immediately from the theorem. Now suppose that $\by \notin \Gamma \bx$ and let $\rho <1/\psi_0$ be given. Let $\psi =1/\rho$, so $\psi>\psi_0$. By the first part, there exists a constant $C\geq 1$ such that for every $\epsilon \in (0,1)$ there exists $\gamma \in \Gamma$ satisfying $\norm \gamma \bx -\by\norm <\epsilon $ and $\norm \gamma \norm <C\epsilon^{-\psi}$; the latter condition implies that  $\epsilon < C^{1/\psi}\norm \gamma \norm^{1/\psi} \leq C^{(n-p)/p}\norm \gamma\norm^{-\rho}$, and hence $\norm \gamma \bx -\by\norm <C^{(n-p)/p}\norm \gamma \norm^{-\rho}$. Since $\by \notin \Gamma \bx$, it follows that the set of $\gamma$  obtained in this way (even corresponding to a sequence of $\epsilon$'s tending to $0$) contains infinitely many distinct elements. This proves the Corollary. \qed 

\medskip
\noindent {\it Proof of Corollary~\ref{newcor}}: The corollary follows immediately from Corollary~\ref{cor:main} and Proposition~\ref{generic}. \qed

\begin{remark}\label{comp-remark}
{\rm In the case $n=2$ and $p=1$, namely the $\Gamma$-action on $\R^2$, Corollary~\ref{cor:main} holds for $\bx$ for which $h(\bx)<2$. We  recall that the result in \cite{LN} for the $\SL (2,\Z)$-action is available for all points which are not multiples of rational vectors, without the condition on exponents. Moreover, for $\bx$ for which $h(\bx)\geq\frac 75$ the value of $\psi$ as in the conclusion exceeds $3$, whereas existence of solutions is assured with $\psi=3$ by the result in \cite{LN} for the action of $\SL (2, \Z)$ and hence that of $\Gamma$. Thus for $\bx$ with $h(\bx)\geq \frac 75$, \cite{LN} offers  better results; however the set of $\bx$ for which that happens has measure $0$. }
\end{remark}
 
Extending further the correspondence as above, we now discuss the analogue of badly approximable matrices, and their significance to our main theorem. 

\begin{definition}\label{proj} 
{\rm Let $1\leq p \leq n-1$ and $\bx \in \R^{(n,p)}$. 
We say that the matrix $\bx$ is {\it projectively badly approximable} if 
 there exists a constant $c(\bx)>0$ such that
$
\norm \omega\bx \norm^p \norm w \norm^{n-p}  > c(\bx)
$ for every $\omega \in \mathcal{M}(n,\Z)\backslash \{0\}$.}
\end{definition}

Badly approximable vectors have been a subject of much study. It would  be worth recalling here the following theorem (cf. \cite{Schmidt2}); see also the note at the end of the section.

\begin{thm}\label{bad}
For $n,p\ge1$,  the set of badly approximable vectors in $\R^{(n,p)}$
is a set of Lebesgue null measure,  of Hausdorff dimension $np$.
\end{thm}

\begin{prop}\label{badly}
Let $1\leq p \leq n-1$ and $q=n-p$. Let $\xi \in \R^{(q, p)}$  and  $\theta \in \GL(p,\R)$. Then the $n\times p$ matrix $\left(\begin{matrix}
\xi \theta \cr
\theta 
\end{matrix} \right)$
is projectively badly approximable if and only if $\xi $ is badly approximable.  
\end{prop}
\proof We shall follow the pattern of the proof of 
Proposition~\ref{corresp}. As in that proposition it suffice to prove the assertion here when $\theta =I$, the identity matrix, as we shall now assume. We shall follow the notation as in Proposition~\ref{corresp}. We note that $\bx$ is projectively badly approximable if and only if 
$$\inf_{\omega \in \mathcal {M} (n,\Z)\setminus \{0\}}\norm \omega \bx\norm\norm \omega \norm^{(n-p)/p}>0,$$  whereas 
$\xi$ is badly approximable if and only if  
 $$\inf_{\alpha \neq 0} d(\alpha)\norm \norm \alpha \norm^{(n-p)/p}>0.$$ The desired assertion therefore follows from the inequalities~(4.1) as in the proof of Proposition~\ref{corresp}, for the value $b=(n-p)/p$. \qed

\medskip

\noindent {\bf Note}: W.M. Schmidt proved (see \cite{Schmidt2}),  apart from Theorem~\ref{bad} as above, stronger results about the class of  badly approximable systems of vectors, in various respects. It should be evident to the interested reader that via the connection described in Proposition~\ref{badly}, correspondingly stronger results could be deduced for projectively badly approximable systems as introduced above. We shall however not go into the details of this here.

\section{ Exponent of diophantine approximation }

For $\bx, \by \in \R^{(n,p)}$, where $1\leq p \leq n-1$,  following \cite{BL} and \cite{LN} we define  the {\it exponent of approximation} of the action of $\Gamma$, corresponding to the pair $(\bx, \by)$, as
$$
e(\bx,\by)=\sup \left\{ \mu \in \R \mid  \norm \gamma \bx - \by \norm < \frac{ 1}{ \norm\gamma \norm^{\mu}} \mbox{ for infinitely many } \gamma \in \Gamma\right\}.
$$ 

In this section we prove the following result, which is  a restatement of Theorem~\ref{thm:exp} stated in the introduction. 

\begin{thm}\label{lowerbound}
Let $1\leq p\leq n-1$.
Then, for Lebesgue almost every pair $(\bx,\by)\in \R^{(n,p)} \times \R^{(n,p)}$, $e(\bx, \by) =(n-p)/ p$.
\end{thm}

\proof As the set of pairs $(\bx, \by)$ such that $\by \notin \Gamma \bx$ is a set of full Lebesgue measure in $\R^{(n,p)} \times \R^{(n,p)}$, it follows immediately from Corollary~\ref{newcor} that $e(\bx, \by)\geq (n-p)/p$ for almost all $(\bx, \by)$. 

Let $ \bx=(x_1, \dots, x_p) $ with $x_1, \dots, x_p$ linearly independent vectors in $\R^n$.
We shall show that $e(\bx, \by)\leq (n-p)/p$ for almost all $\by$.
 The proof of this is along the lines of the proof of the upper bound of the generic density approximation exponent of the linear action of the modular group $\SL(2,\Z)$ on $\R^2$ given in \cite{LN2}, Section~5.

For $\bz \in \R^{(n,p)}$ and $r>0$, let  $B(\bz,r)=\{\by\in \R^{(n,p)}: \norm  \bz - \by \norm <r\} $. It suffices to  show that for any $\rho >0$, $e(\bx, \by)\leq (n-p)/p$ for almost all $\by \in B(0,\rho)$.  Let $\rho >0$ be given and 
$B=B(0,\rho+1)$. Clearly, for $\by \in B(0,\rho)$ and $\mu \geq 0$, if  $\gamma \in \Gamma$ is such that $\norm \gamma \bx -\by \norm <\norm \gamma \norm^{-\mu}$ then $\gamma \bx \in B$. 

We note that 
there exists a positive constant $C>0$ such that for all  $q\in \N$,  
$$
\# \{ \gamma  \in \mathcal{M}(n,\Z) \mid  \norm\gamma \norm\le q \hbox{ \rm and }   \gamma \bx\in B \}  \le C q^{n(n-p)}. \leqno (5.1)
$$
This follows from Minkowski's theorem, since for each $q$, the set as above consists of lattice points in $\{\omega \in  \mathcal{M}(n,\R) \mid  \norm\omega  \norm\le q  \hbox{ \rm and }  \omega \bx \in B\}$, which is a convex symmetric body in the vector space $\mathcal{M}(n,\R)$ whose Lebesgue measure is  $C
 q^{n(n-p)}$, for a suitable constant~$C$. 

Now let $\mu >(n-p)/p$ be given, say $\mu =\frac{n-p}p (1+\delta)$, where $\delta>0$. Let $\ell$ be the standard Lebesgue measure on $\mathcal{M}(n,\R)$. We note that  for any $\bz \in \R^{(n,p)}$ and $r>0$ we have $\ell (B(\bz,r))=2^{np}r^{np}$. 

For $k \ge 1$, let 
$\Gamma_k = \{\gamma \in \Gamma \mid  \norm \gamma \norm = k, \gamma \bx \in B\}$ and $N_k=\# \Gamma_k$, the cardinality of $\Gamma_k$. 
By (5.1) we have $$\displaystyle  N_1+\ldots + N_q \le C q^{n(n-p)} \hbox{ \rm for all } q\in \N.$$
For all $q\in \N$ we have
$$
\ell \left(\bigcup_{k = 1}^q   \bigcup_{\gamma \in \Gamma_k} B(\gamma \bx, k^{-\mu}) \right)
\le 2^{ np} \Sigma_{k=1}^q {N_k \over (k^\mu)^{ np}}=  2^{np}\Sigma_{k=1}^q {N_k \over k^{\mu np}}
$$
$$
=2^{ np} \left( \Sigma_{k=1}^{q-1} (N_1+\ldots N_k) \left({1 \over k^{\mu np}}- {1 \over (k+1)^{\mu np}}\right)
+ {N_1+\ldots +N_q \over q^{\mu np}} \right)
$$
$$
\le  2^{np} C\left( \Sigma_{k=1}^{q-1}   k^{n(n-p)} \left({1 \over k^{\mu np}}- {1 \over (k+1)^{\mu np}} \right)
+ {q^{n(n-p)} \over q^{\mu np}} \right).
$$
Using that
$\displaystyle {1 \over k^{\mu np}}- {1 \over (k+1)^{\mu np}}< \mu np {1 \over k^{1+\mu np}}= \mu np \frac 1{k^{1+n(n-p)\delta}}$, and $\mu = \frac{n-p}p(1+\delta)$ we now obtain 
$$
\ell \left(\bigcup_{k = 1}^q   \bigcup_{\gamma \in \Gamma_k} B(\gamma \bx, k^{-\mu}) \right)
\le  2^{np}C \left( \mu np \Sigma_{k=1}^{q-1}    {1   \over k^{1+n(n-p)\delta}  }
+  {1 \over q^{n(n-p)\delta} }  \right).
$$
As $\displaystyle n(n-p)\delta >0$, it follows that the right hand side term of the above inequality converges as $q \rightarrow \infty$. Thus
$\displaystyle \ell \left(\bigcup_{k \ge 1}   \bigcup_{\gamma \in \Gamma_k} B(\gamma \bx, k^{-\mu}) \right) < \infty$.
Applying the Borel-Cantelli Lemma, we get  that the set 
$$
\limsup_{q \rightarrow \infty}  \bigcup_{\gamma \in \Gamma_q} B(\gamma \bx, q^{-\mu})  
=\bigcap_{q \ge 1} \,  \bigcup_{\gamma \in \Gamma_k, \, k\ge q} B(\gamma \bx, k^{-\mu}) \subset B
$$
 is a null measure set. For any $\by $ in $B(0,\rho)$ which is in the complement of this subset there are only be finitely many $\gamma \in \Gamma $ such that $\by \in B(\gamma \bx , \norm \gamma \norm ^{-\mu})$, namely such that   
 $\norm \gamma \bx -\by\norm <\norm \gamma \norm ^{-\mu}$, and hence $e (\bx, \by) \le \mu$. As this holds for all $\mu >(n-p)/p$ we get that for any $\bx $ as above, $e (\bx, \by) \le (n-p)/p $ for almost all $\by \in B(0,\rho)$. Since this holds for all $\rho >0$ 
 this proves the assertion in the theorem.  \qed
 
 \noindent
{\bf Acknowledgements.} 
We graciously acknowledge the support of R\'egion Provence-Alpes-C\^ote d'Azur through the project APEX {\it Syst\`emes dynamiques: Probabilit\'es et Approximation Diophantienne} PAD, CEFIPRA through the project No. 5801-B and the projet MATHAMSUD No. 38889 DCS: Dynamics of Cantor Systems.

\vskip1cm

\begin{tabular}{ll}
S.G. Dani & Arnaldo Nogueira \\
UM-DAE Centre for Excellence  & Aix Marseille Universit\'e, CNRS, \\
~~~~~in Basic Sciences, & Centrale Marseille,\\
University of Mumbai, Mumbai, India\hspace{15mm} &  I2M, Marseille, France \\[1mm]
shrigodani@cbs.ac.in & arnaldo.nogueira@univ-amu.fr \\

\end{tabular}

\begin{thebibliography}{HD}

\bibitem{BL}
Yann Bugeaud and Michel Laurent.{\it On exponents of homogeneous and inhomogeneous diophantine approximation}, 
Moscow Math. Journal, Vol. 5, Number 4 (2005), 747- 766.



\bibitem{JSD} 
J. S. Dani. Density properties of orbits under discrete groups,
J. Indian Math. Soc. (N.S.) 39 (1975), 189 - 217 (1976). 

\bibitem{DR}
S.G. Dani and S. Raghavan. 
Orbits of Euclidean frames under discrete linear groups, Israel J. Math. 
36 (1980), 300 - 320.

\bibitem{GG}
Anish Ghosh, Alexander Gorodnik and Amos Nevo. Best possible rates of distribution of dense lattice orbits in homogeneous spaces, J. Reine angew. Math. 745 (2018), 155 - 188. 


\bibitem{Ka}
Yonatan R. Katznelson.  Integral matrices of fixed rank, Proc. of the American Mathematical Soc. Volume 120, Number 3 (1994), 667-675.



\bibitem{Ke}
Dubi Kelmer. Approximation of points in the plane by generic lattice orbits, J. Modern Dynamics 11 (2017), 143 - 153.



\bibitem{LN}
Michel Laurent and  Arnaldo Nogueira. 
Approximation to points in the plane by $\SL(2,\Z)$-orbits,
J. Lond. Math. Soc. (2) 85 (2012), no. 2, 409 - 429. 

\bibitem{LN2}
Michel Laurent and  Arnaldo Nogueira. Inhomogeneous approximation with coprime integers and lattice orbits, Acta Arithmetica, 154.4 (2012), 413-427.

\bibitem{MW}
Francois Maucourant and Barak Weiss. Lattice actions on the plane revisited, Geom. Dedicata 157 (2012), 1 - 21. 


\bibitem{meyer} 
Yves Meyer. 
Algebraic numbers and harmonic analysis,
North-Holland Mathematical Library, Vol. 2. North-Holland Publishing Co., Amsterdam-London; American Elsevier Publishing Co., Inc., New York, 1972. x+274~pp. 


\bibitem{Schmidt2} 
Wolfgang Schmidt. {\it Diophantine Approximation}, Lecture Notes in Mathematics 785, Springer-Verlag, 1980.

\end{thebibliography}
\end{document}